\documentclass[10pt,journal]{IEEEtran}
\usepackage{amsmath,graphicx,color,cite,wrapfig,xcolor,bm,svg}
\usepackage{booktabs}
\usepackage{multirow}
\usepackage{graphicx}

\newcommand{\beq}{\begin{equation}}
\newcommand{\eeq}{\end{equation}}

\newcommand{\bdelta}{{\bm{\delta}}}

\newcommand{\bftheta}{{\mbox{\boldmath $\theta$}}}

\newcommand{\cG}{{\mbox{$\mathcal{G}$}}}

\newtheorem{example}{Example}

\def\remark{\addtocounter{remark}{1}\def\@currentlabel{\theremark}
\emph{Remark~\theremark}. } \makeatother
\newcounter{remark}

\newcommand{\bx}{\mathbf{x}}
\newcommand{\bw}{\mathbf{w}}

\newcommand{\bX}{\mathbf{X}}
\newcommand{\bS}{\mathbf{S}}

\newcommand{\bZ}{\mathbf{Z}}

\newcommand{\bA}{\mathbf{A}}
\newcommand{\bB}{\mathbf{B}}

\newcommand{\bz}{\mathbf{z}}
\newcommand{\bu}{\mathbf{u}}
\newcommand{\bv}{\mathbf{v}}
\newcommand{\by}{\mathbf{y}}

\usepackage{amsmath}
\usepackage{amssymb}
\usepackage{flushend}
\usepackage{cancel}
\usepackage[normalem]{ulem}
\usepackage{multirow}
\usepackage{caption}
\usepackage{makecell}
\usepackage{threeparttable}

\usepackage[margin=1in]{geometry}
\title{Non-convex Min-Max Optimization: Applications, Challenges, and Recent Theoretical Advances}

\author{{Meisam Razaviyayn, Tianjian Huang, Songtao Lu, Maher Nouiehed, Maziar Sanjabi, Mingyi Hong}
}
	
\begin{document}

\maketitle

\begin{abstract}
The min-max optimization problem, also known as the saddle point problem, is a classical optimization problem which is also studied in the context of zero-sum games. Given a class of objective functions, the goal is to find a value for the argument  which leads to a small objective value even for the worst-case function in the given class. Min-max optimization problems have recently become very popular in a wide range of signal and data processing applications such as fair beamforming, training  generative adversarial networks (GANs), and robust machine learning (ML), to just name a few. The overarching goal of this article is to provide a survey of recent advances for an important subclass of min-max problem in which the minimization and maximization problems can be non-convex and/or non-concave. 
In particular, we first present a number of applications to showcase the importance of such min-max problems; then, we discuss key theoretical challenges, and provide a selective review of some exciting recent theoretical and algorithmic advances in tackling non-convex min-max problems. Finally, we point out open questions and future research directions.\footnote{The manuscript is accepted in IEEE Signal Processing Magazine.} 
\end{abstract}

{\section{Introduction}
\label{sec:intro} 
Recently, the class of  non-convex min-max optimization problems has attracted significant attention across signal processing, optimization, and ML communities. The overarching goal of this paper is to provide a selective survey of the applications of such a new class of problem, discuss theoretical and algorithmic challenges, and present some recent advances in various directions. 

To begin our discussion, let us consider the following generic problem formulation: 
\begin{align}\label{eq:MinMax}
& \min_{{\bx}}\max_{{\by}} \quad f(\bx,\by) \tag{Min-Max}\\
& \; \mbox{s.t.} \quad \bx\in \mathcal{X}\subseteq\mathbb{R}^d, \;\; \by\in \mathcal{Y}\subseteq\mathbb{R}^b,\nonumber
\end{align}
where $f(\cdot,\cdot):\mathbb{R}^{d}\times\mathbb{R}^{b}\rightarrow \mathbb{R}$ is differentiable with Lipschitz continuous gradient in $(\bx,\by)$, possibly non-convex in $\bx$ and possibly non-concave in $\by$; 
$\bx\in\mathbb{R}^{d}$ and $\by\in\mathbb{R}^{b}$ are the optimization variables; $\mathcal{X}$ and $\mathcal{Y}$ are the feasible sets, which are assumed to be closed and convex. Notice that, while we present this article around the above min-max formulation, extending the ideas and discussions to max-min problems is straight forward. 

When problem~\eqref{eq:MinMax} is convex in $\bx$ and concave in $\by$, the corresponding variational inequality (VI) becomes monotone, and a wide range of algorithms have been proposed for solving this problem; see, e.g., \cite{nemirovski2004prox,juditsky2016solving,  tseng2008accelerated, nesterov2007dual}, and the references therein. 
However, as we will discuss in this article, solving  min-max problems is challenging in non-convex setting. Such non-convex min-max optimization problems appear in  different applications  in signal processing   (e.g., robust transceiver design,
fair resource allocation \cite{liu11ICC}, 
communication in the presence of jammers {\cite{gohary09}}), distributed signal processing~\cite{Giannakis15,nedic18proceeding}), 
and ML (e.g., 
robust training of neural networks~\cite{madry2018}, 
training generative adversarial networks (GANs)~\cite{goodfellow2014generative, arjovsky2017wasserstein}, 
and fair inference~\cite{zhang2018mitigating,baharlouei2019r}. 
More generally, any design problem in the presence of model uncertainty or adversary can be modeled as an optimization of the form~\eqref{eq:MinMax}. In this setup, $\bx$ is the design parameter that should be optimized, while $\by$ is the uncertainty/adversary parameter which is not accurately measured, or may be adjusted by an adversary. In such scenarios, the goal in formulation~\eqref{eq:MinMax} is to find a solution~$\bx = \bar{\bx}$ that has a robust performance against all uncertainty/adversary values of $\by\in\mathcal{Y}$. Such a robustness requirement has long been deemed important in signal processing community, and it has recently played a crucial role in designing modern ML tools. 

Despite the rising interests in   nonconvex min-max problems, seldom have they been rigorously analyzed in either classical optimization or signal processing literature.  In this article, we
first present a number of applications to showcase the importance of such min-max problems, then we discuss key theoretical challenges, and provide a selective review of some recent theoretical and algorithmic advances in tackling the class of non-convex min-max problems. Finally, we will point out open questions and future research directions. 

\section{Applications of non-convex min-max problems}
\label{sec:applications} 

To appreciate the importance of problem~\eqref{eq:MinMax}, let us first present a number of key applications of non-convex min-max optimization problems. 
\medskip
\subsubsection{Generative adversarial networks (GANs)}
GANs~\cite{goodfellow2014generative} have recently gained tremendous popularity due to their unique ability to learn complex distributions and generate realistic samples, e.g., high resolution fake  images.  
In the absence of labels,  GANs  aim at finding a mapping from a known distribution, e.g. Gaussian, to an unknown data distribution, which is only represented by empirical samples~\cite{goodfellow2014generative}.


GANs consist of two neural networks: the generator and the discriminator. The goal of the generator is to generate \textit{fake} samples which look like \textit{real} samples in the distribution of interest. This process is done by taking i.i.d. samples from a known distribution such as Gaussian and transform it to samples similar to real ones via trained neural network. On the other hand, the discriminator's objective is to correctly classify the fake samples generated by the generator and the real samples drawn from the distribution of interest. The  two-player game between the generator and the discriminator can be modeled as a min-max optimization problem \cite{goodfellow2014generative}:
\begin{equation}
\label{eq:GAN}
\min_{\bw_g}  \; \max_{\bw_d} \;\;V(\bw_g,\bw_d),
\end{equation}
where $\bw_g$ is the generator's parameter; $\bw_d$ is the discriminator's parameter; and $V(\cdot,\cdot)$ shows the cost function of the generator (which is equal to the negative of the discriminator's cost function). 
This min-max objective can be also justified as minimizing some distance between the distribution of real samples and the distribution of generated samples.  In this interpretation, the distance between the two distributions is computed by solving a maximization (dual) problem; and the goal is to minimize the distance between the distribution of generated samples and the distribution of real samples.  Various distance measures have been used for training GANs such as Jensen-Shannon divergence \cite{goodfellow2014generative}, $f$-divergence, 
and Wasserstein distance \cite{arjovsky2017wasserstein}. All these distances lead to non-convex non-concave min-max formulations for training GANs.  

\medskip
\subsubsection{Fair ML}
The past few years have witnessed several reported instances of ML algorithms suffering from systematic discrimination against individuals of certain protected groups; see, e.g., \cite{baharlouei2019r,xu2018fairgan,DBLP:conf/icml/MadrasCPZ18}, and the references therein. Such instances stimulate strong interest in the field of fairness in ML which in addition to the typical goal of having an accurate learning model, brings fairness to the learning task. Imposing fairness on ML models can be done through three main approaches: preprocessing approaches, in-processing approaches, and postprocessing approaches.

To understand these three approaches, consider a ML task over a given random variables $\bX \in \mathbb{R}^{d}$ representing the non-sensitive data attributes and $\bS \in \mathbb{R}^{k}$ representing the sensitive attributes (such as age, gender, ethnicity, etc.). Pre-processing approaches tend to hinder discrimination by masking the training data before passing it to the decision making process. Among these methods,  recent works~\cite{xu2018fairgan,DBLP:conf/icml/MadrasCPZ18} have used an adversarial approach which seeks to learn a data representation $\bZ = \zeta(\bX, \bS)$ capable of minimizing the loss over the classifier~$g( \bZ)$, and protecting the sensitive attributes $\bS$ from an adversary $h(\bZ)$ that tries to reconstruct $\bS$ from $\bZ$. This requires solving the following min-max optimization problem: 
\[
\min_{\zeta, g} \max_{h} \, \mathbb{E}_{\bX, \bS} \{{\cal L}(\zeta, g, h)\}.
\]
Realizing the functions as neural networks, this formulation leads to non-convex min-max optimization problem. 

Contrary to pre-processing methods, in-processing approaches  impose fairness during training procedure. For example, they impose fairness by adding a regularization term that penalizes statistical dependence between the learning model output and the sensitive attributes~$\bS$. Let $g_{\bftheta}(\bX, \bS)$ be a certain output of the learning model. One can balance the learning accuracy and fairness by solving the following optimization problem:

\begin{equation}\label{Fair_Machine_Learning_Formulation}
\min_{\bftheta} \;\;\mathbb{E}\left\{L\left(\bftheta, \bX \right) \right\} + \lambda \; \rho\left(g_{\bftheta}(\bX, \bS), \bS\right),
\end{equation}
where $\rho(\cdot, \cdot)$ is a statistical independence measure and $L(\cdot,\cdot)$ denotes the training loss function. For example, in the classification task in which $\bX$ contains both the input feature and the target variable, the function~$L(\cdot,\cdot)$ measures the classification error of the trained classifier. Here, the parameter $\lambda$ is a positive scalar balancing fairness and accuracy of the output model. When $\lambda \rightarrow \infty$, this optimization problem focuses more on making $g_{\bftheta}(\bX, \bS)$ and $\bS$ independent, resulting in a fair inference. However, when $\lambda = 0$, no fairness is imposed and the focus is to maximize the accuracy of the model output.

Various statistical dependence measures have been proposed for use in this formulation.  For example, \cite{baharlouei2019r} proposed using R\'{e}nyi correlation to impose fairness. The R\'{e}nyi correlation between two random variables $\bA$ and $\bB$ is defined as $\rho(\bA,\bB) \triangleq \sup_{k,\ell} \rho_p(k(\bA),\ell(\bB))$ where $\rho_p$ is the Pearson correlation coefficient and
the supremum is over the set of measurable functions $k(\cdot)$ and $\ell(\cdot)$. Plugging the definition of R\'{e}nyi correlation in~\eqref{Fair_Machine_Learning_Formulation} leads to a natural min-max formulation, which is the focus of this article.

\medskip
\subsubsection{Adversarial ML}\label{sub:attack}
The formulation \eqref{eq:MinMax} is also instrumental to model the dynamic process of {\it adversarial learning}, where the model training process involves some kind of ``adversary". Depending on whether the goal is to  break the ML model or to make it more robust, one can formulate different min-max optimization problems, as we briefly discuss in the following sections.

\smallskip
\noindent{\bf Adversarial attacks.} First, let us take the viewpoint of the adversary, who would like to break a ML model so that it is more likely to produce wrong predictions. In this scenario, the adversary tries to increase the error of a well-trained ML model; therefore its behavior is modeled as the {\it outer} optimization problem, aiming to reduce the performance of the trained model. On the other hand, the training process is modeled as the {\it inner} optimization problem aiming to minimize the training error. 

To be more specific, take the poisoning attack~\cite{zominmax19} as an example. Let $\mathcal D := \{ \mathbf u_i, t_i \}_{i=1}^N $ denote   the training dataset, where $\mathbf u_i$ and $t_i$ represent the features and target labels of sample $i$ respectively. Each data sample~$\mathbf{u}_i$ 
can be corrupted by a perturbation vector $\boldsymbol{\delta}_i$ to generate a ``poisoned" sample $\mathbf u_i + \boldsymbol{\delta}_i$.  
Let $\boldsymbol{\delta}:= (\boldsymbol{\delta}_1,\ldots,\boldsymbol{\delta}_N)$ be the collection of all poisoning attacks. Then, the poisoning attack problem is  formulated as 
{
\begin{align}\label{eq: poison_robust_attack}
\displaystyle \max_{\boldsymbol{\delta}\text{:}\;\|\boldsymbol{\delta}_i\|\leq \varepsilon} \min_{\bw}  \; 
\sum_{i=1}^{N}{\ell}(p(\bu_i+\boldsymbol{\delta}_i;\bw),t_i) 
\end{align}
}%
where $\bw$ is the weight of the neural network; 
$p(\cdot)$ is the predicted output of the neural network; and $\ell(\cdot)$ is the loss function.     
The constraint $\| \boldsymbol{\delta}_i \| \leq \varepsilon$ indicates that the poisoned samples should not be too different from the original ones, so that the attack is not easily  detectable. Note that the ``max-min'' problem \eqref{eq: poison_robust_attack} can be written equivalently in  the form of \eqref{eq:MinMax} by adding a minus sign to the objective. 

\smallskip
\noindent{\bf Defense against adversarial attacks.} It has been widely observed that ML models, especially neural networks, are highly vulnerable to adversarial attacks, including the poisoning attack discussed in the previous subseciton, or other popular attackes such as Fast Gradient Sign Method (FGSM) attack~\cite{FGSM} and Projected Gradient Descent (PGD) attack~\cite{PGD}. These adversarial attacks show that a small perturbation in the data input can significantly change the output of a neural network and deceive different neural network architectures in a wide range of applications. 
To make ML models robust against adversarial attacks, one popular approach is to solve the following 
robust training problem \cite{madry2018} (using similar notations as in \eqref{eq: poison_robust_attack}):
\begin{align} \nonumber
\min_{\bw }\, \; \sum_{i=1}^{N}\;\max_{\boldsymbol{\delta}\text{:}\;\|\boldsymbol{\delta}_i\|\leq \varepsilon}   {\ell}(p(\bu_i+\boldsymbol{\delta}_i;\bw ), t_i).
\end{align}
Note that compared with \eqref{eq: poison_robust_attack}, the roles of minimization and maximization have been switched.
Clearly, this optimization problem is of the form~\eqref{eq:MinMax}.

\medskip
\subsubsection{Distributed processing}
Some constrained non-convex optimization problems could also be formulated as a min-max saddle point problem by leveraging the primal dual approach or the method of Lagrange multipliers.  An example of that appears in distributed data processing over networks.
Consider a network of $N$ nodes in a graph ${\mathcal{G}}=\{{\mathcal{V}}, {\mathcal{E}}\}$ with $|{\mathcal{V}}|=N$ vertices. The nodes can communicate with their neighbors, and their goal is to jointly solve the optimization problem:
\[
\min_{z} \quad \sum_{i=1}^{N}  g_{i}(z),
\]
where each $g_{i}(\cdot)$ is a smooth function only known by node $i$. Further, for simplicity of presentation, assume that $z \in \mathbb{R}$.

Such a distributed optimization setting has been widely studied in the optimization and signal processing communities over the past few decades. Let $x_{i}$ be node $i$'s local copy of $z$. A standard first step in distributed  optimization is to rewrite the above problem as:
{
	\begin{align}\label{eq:distributed}
	\min_{{x} \in \mathbb{R}^{N}} g({\bx}):= \sum_{i=1}^{N} g_{i}(x_{i}) \quad \quad \textrm{ s.t.} \quad \mathbf{A}{\bx}=0,
	\end{align}
	}%
where ${\bA} \in \mathbb{R}^{|\mathcal{E}| \times N}$ is the  incidence matrix for graph $\cG$ and $\bx = (x_1,\ldots,x_N)$ is the concatenation of all copies of the decision variable. The linear constraint in \eqref{eq:distributed} enforces $x_i=x_j,$ if $i,j$ are neighbors. 
Problem \eqref{eq:distributed} can be rewritten as\footnote{It can be shown that finding a stationary solution of \eqref{eq:reform} is equivalent to finding a stationary solution for \eqref{eq:distributed}; see \cite{minmaxtsp19}.}
{
	\begin{align}\label{eq:reform}
	\max_{{y}\in\mathbb{R}^{|\mathcal{E}|}} \min_{{x}\in\mathbb{R}^{N}}\sum_{i=1}^{N}g_i(x_i)+ {\by}^T{\bA}{\bx}
	\end{align}}%
	where ${\by}$ is the Lagrangian multiplier. Clearly, \eqref{eq:reform} is in the form of \eqref{eq:MinMax}, where the coupling between $\bx$ and $\by$ is {\it linear}. A number of algorithms have been developed for it; 
 see a recent survey \cite{Chang20SPM}.

\medskip
\subsubsection{Max-Min fair transceiver design}
Consider the problem of resource allocation in a wireless communication system, where  $N$~transmitter-receiver pairs are communicating. The goal is to maximize the minimum rate among all users. To be specific, consider a setting with $K$ parallel channels. User $i$ transmits with power $\mathbf{p}_i:=[p^1_i; \cdots, p^K_i]$, and its rate is given by: $r_i(\mathbf{p}_1, \ldots, \mathbf{p}_N) = \sum_{k=1}^{K}\log\bigg(1 + \dfrac{a^k_{ii} p^k_i}{\sigma^2_i + \sum_{j=1, j \neq i}^{N} a^k_{ji} p^k_j}\bigg)$ (assuming  Gaussian signaling), which is a non-convex function in $\mathbf{p}$. Here  $a^k_{ji}$ denotes the channel gain between transmitter $j$ and receiver $i$ on the $k$-th channel, and $\sigma^2_i$  is the noise power of user~$i$. Let $\mathbf{p}:=[\mathbf{p}_1;\cdots; \mathbf{p}_N]$, then the max-min fair power control problem is given by \cite{liu11ICC}
\begin{equation} \label{eq:max-min-fair-discrete}
    \max_{\mathbf{p} \in \mathcal{P}} \;\min\limits_{i} \; \{r_{i}(\mathbf{p})\}_{i=1}^N,
\end{equation}
where  $\mathcal{P}$ 
denotes the set of feasible power allocations. While the inside minimization is over a discrete variable~$i$, we can reformulate it as a minimization over continuous variables using transformation:
\begin{equation}\label{eq:max-min-fair}
\max\limits_{\mathbf{p}\in\mathcal{P}} \;\min\limits_{\by \in \Delta}\quad  \sum\limits_{i=1}^{N} r_{i}(\mathbf{p}_1,\cdots, \mathbf{p}_N) \times y_i,
\end{equation}
where $\Delta \triangleq \{\by\,|\, \by \geq \mathbf{0};\; \sum_{i=1}^N y_i = 1 \}$ is the probability simplex. Notice that the inside minimization problem is linear in $\by$. Hence, there always \emph{exists} a solution at one of the extreme points of the simplex $\Delta$. Thus, the formulation~\eqref{eq:max-min-fair} is equivalent to the formulation~\eqref{eq:max-min-fair-discrete}. By multiplying the objective by the negative sign, we can transform the above ``max-min" formulation to ``min-max" form  consistent with \eqref{eq:MinMax}, i.e., 
\begin{equation}\nonumber
\min\limits_{\mathbf{p}\in\mathcal{P}}\; \max\limits_{\by \in \Delta}\quad  \sum\limits_{i=1}^{N} -r_{i}(\mathbf{p}_1,\cdots, \mathbf{p}_N) \times y_i 
\end{equation}


\medskip
\subsubsection{Communication in the presence of jammers}\label{sub:jammer}
Consider a variation of the above problem, where $M$ jammers participate in an $N$-user $K$-channel interference channel transmission. The jammers' objective is to reduce the sum rate of the system by transmitting noises, while the goal for the regular users is to transmit as much information as possible. 
We use $p^k_i$ (resp. $q^k_j$) to denote the $i$-th regular user's (resp. $j$-th jammer's) power on the $k$-th channel. The corresponding sum-rate max-min problem can be formulated as: 
	\begin{align}\label{eq:jammer}
	& \max_{\mathbf{p}}\;\min_{\mathbf{q}}\;
	 \sum_{k,i,j}  \log \left( 1 + \dfrac{a^k_{ii} p_{i}^{k}}{\sigma_{i}^{2} + \sum_{\ell=1,j \neq i}^{{N}} a_{\ell i}^{k} p_{\ell}^{k} + b_{ji}^{k}q_{j}^{k}} \right), \nonumber\\
	 & \;\mbox{s.t.}  \;\; \mathbf{p} \in \mathcal{P},\; \; {\mathbf{q}\in \mathcal{Q}, }
	\end{align}
where $a_{\ell i}^k$  and $b_{ji}^k$ represent the $k$-th channels between the regular user pairs~$(\ell, i)$ and regular and jammer pair~$(i,j)$, respectively.  
Here $\mathcal{P}$ and $\mathcal{Q}$ denote the set of feasible power allocation constraints for the users and the jammers. Many other related formulations have been considered, mostly from the game theory perspective~\cite{gohary09}. Similar to the previous example, by multiplying the objective by a negative sign, we obtain an optimization problem of the form~\eqref{eq:MinMax}. 

}

\section{Challenges}\label{sec:challenges}

Solving min-max  problems even up to simple notions of stationarity could be extremely challenging in the non-convex setting. This is not only because of the non-convexity of the objective (which prevents us from finding global optima), but also is due to aiming for finding a min-max solution.
To see the challenges of solving  non-convex min-max problems, let us compare and contrast  the optimization problem~\eqref{eq:MinMax} with the regular smooth non-convex optimization problem:
\begin{equation}\label{eq:GeneralNon-convex}
\min_{\bz \in \mathcal{Z}} \;\; h(\bz). 
\end{equation}
 where the gradient of function $h$ is Lipschitz continuous.

While solving general non-convex optimization problem~\eqref{eq:GeneralNon-convex} to global optimality is hard, one can apply simple iterative algorithms such as projected gradient descent  (PGD) to~\eqref{eq:GeneralNon-convex} by running the iterates
\[
\bz^{r+1} = \mathcal{P}_{\mathcal{Z}}(\bz^r - \alpha \nabla h(\bz^r)),
\]
where $r$ is the iteration count; $\mathcal{P}_{\mathcal{Z}}$ is projection to the set $\mathcal{Z}$; and $\alpha$ is the step-size. 
Algorithms like PGD enjoy two properties: 
\begin{itemize}
    \item[i)] The quality of the iterates improve over time, i.e., $h(\bz^{r+1})\leq h(\bz^r)$, where~$r$ is the iteration number. 
    \item[ii)] These algorithms are guaranteed to converge to (first-order) stationary points with global iteration complexity guarantees~\cite{nesterov2013introductory} under a mild set of assumptions.
\end{itemize}   

\begin{figure}
	\begin{center}
		\includegraphics [width=0.5\textwidth]{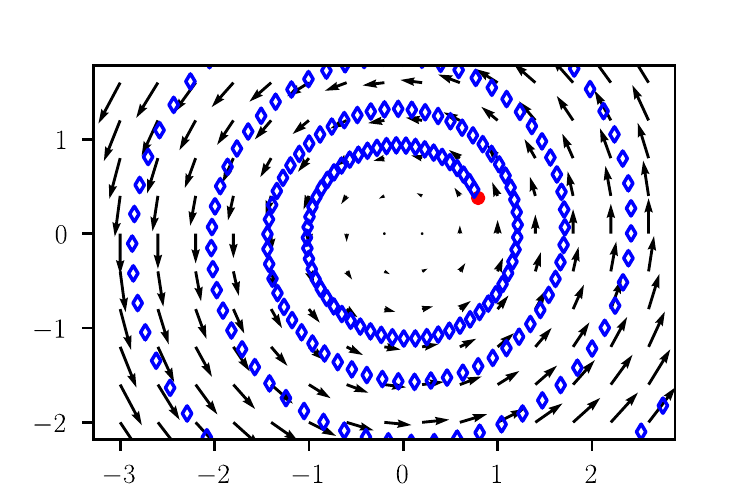}
		\caption{\footnotesize GDA trajectory for the function $f(x,y) = xy$. The iterates of GDA diverge even in this simple scenario. The GDA algorithm starts from the red point and moves away from the origin (which is the optimal solution).}
		\label{Nowazin17}
	\end{center}
\end{figure}

The above two properties give enough confidence to researchers to apply projected gradient descent to many non-convex problems of the form~\eqref{eq:GeneralNon-convex}  and expect to find ``reasonably good'' solutions in practice.  In contrast, \textit{there is no widely accepted optimization tool} for solving general non-convex min-max problem~\eqref{eq:MinMax}. A simple extension of the PGD to the min-max setting is the gradient-descent ascent algorithm (GDA). This popular algorithm simply alternates between a gradient descent step on $\bx$ and a gradient ascent step on~$\by$ through the update rules
\begin{align}
    \bx^{r+1}  &= \mathcal{P}_{\mathcal{X}} (\bx^r - \alpha \nabla_x f(\bx^r,\by^r)), \nonumber\\
     \by^{r+1} &= \mathcal{P}_{\mathcal{Y}}(\by^r + \alpha \nabla_y f(\bx^{r},\by^r)), \nonumber
\end{align}
where $\mathcal{P}_{\mathcal{X}}$ and $\mathcal{P}_{\mathcal{Y}}$ are the projections to the sets $\mathcal{X}$ and $\mathcal{Y}$, respectively. The update rule of $\bx$ and $\by$ can be done alternatively as well [i.e., $\by^{r+1} = \mathcal{P}_{\mathcal{Y}}(\by^r + \alpha \nabla_y f(\bx^{r+1},\by^r))$]. Despite popularity of this algorithm, it fails in many practical instances. Moreover, it is not hard to construct very simple examples for which this algorithm fails to converge to any meaningful point; see Fig.~\ref{Nowazin17} for an illustration. 

\section{Recent developments for solving non-convex min-max problems}
\label{sec:Algorithms}
To understand some of the recently developed algorithms for solving non-convex min-max problems, we first need to review and discuss stationarity and optimality conditions for such problems. Then, we highlight some of the ideas leading to  algorithmic developments.

\subsection{Optimality Conditions}\label{sub:optimality}
Due to the non-convex nature of problem~\eqref{eq:MinMax}, finding the global solution is NP-hard in general~\cite{Murty1987}. Hence, the developed algorithms in the literature aimed at finding \textit{``stationary solutions''} to this optimization problem. One approach for defining such stationarity concepts is to look at problem~\eqref{eq:MinMax} as a game. In particular, one may ignore the order of minimization and maximization in problem~\eqref{eq:MinMax} and view it as a zero-sum game between two players. In this game, one player is interested in solving the problem: 
$\min_{\bx \in \mathcal{X}} \quad f(\bx,\by)$, 
while the other player is interested in solving:
$\max_{\by \in \mathcal{Y}} \; f(\bx,\by).$
Since the objective functions of both players are non-convex in general, finding a global Nash Equilibrium is not computationally tractable~\cite{Murty1987}. Hence, we may settle for finding a point satisfying first-order optimality conditions for each player's objective function, i.e., finding a point $(\bar{\bx},\bar{\by})$ satisfying
\begin{equation}\tag{Game-Stationary} \label{eq:Game-Stat}
\begin{split}
\langle \nabla_{\bx} f(\bar{\bx},\bar{\by}) ,\bx - \bar{\bx} \rangle &\geq 0,\;\forall \bx \in \mathcal{X}\\
&\textrm{and}\\
\langle \nabla_{\by} f(\bar{\bx},\bar{\by}) ,\by - \bar{\by} \rangle &\leq 0,\;\forall \by \in \mathcal{Y}.
\end{split}
\end{equation}
This condition, which is also referred to as ``quasi-Nash Equilibrium" condition in \cite{pang2011nonconvex} or ``First-order Nash Equilibrium" condition in \cite{nouiehed2019solving}, is in fact the solution of the VI corresponding to the min-max game. Moreover, one can use fixed point theorems and show existence of a point satisfying~\eqref{eq:Game-Stat} condition under a mild set of assumptions; see, e.g., \cite[Proposition 2]{pang20164}. In addition to existence, it is always easy to \textit{check} whether a given point satisfies the condition~\eqref{eq:Game-Stat}. The ease of checkability and the game theoretic interpretation of the above~\eqref{eq:Game-Stat} condition have attracted many researchers to focus on developing algorithms for finding a point satisfying this notion; see, e.g., \cite{pang2011nonconvex, nouiehed2019solving, pang20164}, and the references therein.

A potential drawback of the above stationarity notation is its ignorance to the order of the minimization and maximization players. 
Notice that the Sion's min-max theorem shows that when $f(\bx,\by)$ is convex in $\bx$ and concave in~$\by$ the minimization and maximization can interchange in \eqref{eq:MinMax}, under the mild additional assumption that either $\mathcal{X}$ or $\mathcal{Y}$ is compact. However, for the general non-convex problems, the minimization and maximization cannot interchange, i.e., 
\[
\min_{\bx \in \mathcal{X}} \max_{\by \in \mathcal{Y}} \;f(\bx,\by) \;\;\neq \;\;\max_{\by \in \mathcal{Y}} \min_{\bx \in \mathcal{X}} \;f(\bx,\by).
\]
Moreover, the two problems may have different solutions.
Therefore, the \eqref{eq:Game-Stat} notion might not be practical in applications in which the minimization and maximization order is important, such as defense against adversarial attacks to neural networks (as discussed in the previous section). To modify the definition and considering the minimization and maximization order, one can define the stationarity notion by rewriting the optimization problem~\eqref{eq:MinMax} as
\begin{equation}\label{eq:OptViewpoint}
\min_{\bx \in \mathcal{X}}\;\;   g(\bx)
\end{equation}
where $g(\bx) \triangleq \max_{\by \in \mathcal{Y}} f(\bx,\by)$ when $\bx \in \mathcal{X}$ and $g(\bx) = +\infty$ when $\bx \notin \mathcal{X}$. Using this viewpoint, we can define a point $\bar{\bx}$ as a stationary point of \eqref{eq:MinMax} if $\bar{\bx}$ is a first-order stationary point of the non-convex non-smooth optimization~\eqref{eq:OptViewpoint}. In other words,
\begin{equation}\tag{Optimization-Stationary}\label{eq:Opt-Stationary}
\mathbf{0} \in \partial g(\bar{\bx}),
\end{equation}
 where $ \partial g(\bar{\bx})$ is Fr\'{e}chet sub-differential of a function $g(\cdot)$ at the point  $\bar{\bx}$, i.e., $\partial g(\bar{\bx}) \triangleq \{ \mathbf{v} \,|\, \lim\inf_{\bx' \mapsto \bx} \big( g(\bx') - g(\bx) - \langle \mathbf{v}, \bx'-\bx\rangle\big) / \big(\|\bx'-\bx\|\big)
\geq 0\}$. It is again not hard to show existence of  a point satisfying~\eqref{eq:Opt-Stationary} under a mild set of assumptions such as compactness of the feasible set and continuity of the function $f(\cdot,\cdot)$. This is because of the fact that any continuous function on a compact set attains its minimum. Thus, at least the global minimum of the optimization problem~\eqref{eq:OptViewpoint} satisfies the optimality condition~\eqref{eq:Opt-Stationary}.
 This is in contrast to the~\eqref{eq:Game-Stat} notion in which even the global minimum of \eqref{eq:OptViewpoint} may not satisfy \eqref{eq:Game-Stat} condition. The following example, which is borrowed from \cite{jin2019minmax}, illustrates this fact.

\medskip
\begin{example}\label{ex:exampleChiJin}  Consider the optimization problem \eqref{eq:MinMax} where the function $f(x,y) = 0.2xy-\cos(y)$ in the region $[-1,1] \times [-2\pi,2\pi]$. It is not hard to check that  this min-max optimization problem has two global solutions $(x^*,y^*) = (0,-\pi)$ and $(0,\pi)$. However, neither of these two points satisfy the condition \eqref{eq:Game-Stat}. \hfill $\blacksquare$ 
\end{example}

 One criticism of the \eqref{eq:Opt-Stationary} notion is the high computational cost of its evaluation for general non-convex problems. More precisely, unlike the \eqref{eq:Game-Stat} notion, checking \eqref{eq:Opt-Stationary} for a given point $\bar{\bx}$ could be computationally intractable for general non-convex  function $f(\bx,\by)$. Finally, it is worth mentioning that, although the two stationary notions \eqref{eq:Opt-Stationary} and $\eqref{eq:Game-Stat}$ lead to different definitions of stationarity (as illustrated in Example~\ref{ex:exampleChiJin}), the two notions could {coincide} in special cases such as when the function $f(\bx,\by)$ is concave in $\by$ and its gradient is Lipschitz continuous; see \cite{lin2019gradient} for more detailed discussion.  
 
 \subsection{Algorithms Based on Potential Reduction}\label{sub:potential}
Constructing a potential and developing an algorithm to optimize the potential function is a popular way to solve different types of games. A natural potential to minimize is the function $g(\bx)$, defined in the previous section. 
 In order to solve \eqref{eq:OptViewpoint} using standard first-order algorithms, one needs to have access to the (sub-)gradients of the function $g(\cdot)$. While presenting the function $g(\cdot)$ in closed-form may not be possible, calculating its gradient at a given point may still be feasible via Danskin's theorem stated below.

\smallskip
\noindent\textbf{Danskin's Theorem \cite{danskin1966theory}:} Assume the function $f(\bx,\by)$ is differentiable in $\bx$ for all $\bx \in \mathcal{X}$. Furthermore, assume that $f(\bx,\by)$  is strongly concave in $\by$ and that $\mathcal{Y}$ is compact. Then, the function $g(\bx)$ is differentiable in $\bx$. Moreover, for any $
\bx \in \mathcal{X}$, we have
$
\nabla g(\bx) =  \nabla_{\bx} f(\bx,\by^*(\bx)),
$
 where  $\by^*(\bx) = \arg\max_{\by\in \mathcal{Y}}f(\bx,\by)$. 
 

This theorem states that one can compute the gradient of the function $g(\bx)$ through the gradient of the function $f(\cdot,\cdot)$ when the inner problem is strongly concave, i.e. $f(\bx,\cdot)$ is strongly concave for all $\bx$. Therefore, under the assumptions of Danskin's theorem, to apply gradient descent algorithm  to \eqref{eq:OptViewpoint}, one needs to run the following iterative procedure:
\begin{subequations}\label{eq:gradientDescent}
\begin{align}
\by^{r+1} &= \arg\max_{\by \in \mathcal{Y}} \;f(\bx^r,\by) \\
\bx^{r+1} &= \mathcal{P}_{\mathcal{X}} ( \bx^r - \alpha  \nabla f(\bx^r,\by^{r+1})).
\end{align}
\end{subequations}
More precisely, the dynamics obtained in~\eqref{eq:gradientDescent} is equivalent to the gradient descent dynamics $\bx^{r+1} = \mathcal{P}_{\mathcal{X}} ( \bx^r - \alpha \nabla g(\bx^r)),$ according to  Danskin's theorem. 
Notice that computing the value of $\by^{r+1}$ in \eqref{eq:gradientDescent} requires finding the exact solution of the optimization problem in (\ref{eq:gradientDescent}a).  In practice, finding such an exact solution may not be computationally possible. Luckily, even an inexact 
version of this algorithm is guaranteed to converge as long as the point $\by^{r+1}$ is computed accurately enough in (\ref{eq:gradientDescent}a). In particular, \cite{jin2019minmax} showed that the iterative algorithm
\begin{subequations}\label{eq:ApproxgradientDescent}
\begin{align}
&\textrm{Find } \by^{r+1} \textrm{ s.t. } f(\bx^r,\by^{r+1}) \geq \max_{\by \in \mathcal{Y}} \;f(\bx^r,\by) -\epsilon\\
&\bx^{r+1} = \mathcal{P}_{\mathcal{X}}(\bx^r - \alpha  \nabla f(\bx^r,\by^{r+1})).
\end{align}
\end{subequations}
is guaranteed to find an ``approximate stationary point" where the approximation accuracy depends on the value of $\epsilon$. Interestingly, even the strong concavity assumption could be relaxed for convergence of this algorithm as long as step (\ref{eq:ApproxgradientDescent}a) is computationally affordable (see \cite[Theorem 35]{jin2019minmax}). 
The rate of convergence of this algorithm is accelerated for the case in which the function $f(\bx,\by)$ is concave in $\by$ (and general nonconvex in $\bx$) in \cite{nouiehed2019solving} and further improved in \cite{NIPS2019_9430} through a proximal-based acceleration procedure. These works  (locally) create strongly convex approximation of the function by adding proper regularizers, and then apply accelerated iterative first-order procedures for solving the approximation. The case in which the function $f(\bx,\by)$ is concave in $\by$ has also applications in solving finite max problems of the form: 
\begin{equation}\label{eq:finiteMax}
\min_{\bx \in \mathcal{X}} \;\max \;\{f_1(\bx),\ldots,f_n(\bx)\}.
\end{equation}
This is due to the fact that this optimization problem can be rewritten as
\begin{equation}\label{eq:finiteMax2}
\min_{\bx \in \mathcal{X}} \;\max_{\by \in \Delta} \;\sum_{i=1}^n y_if_i(\bx),
\end{equation}
where $\Delta \triangleq \{\by\,|\, \by \geq \mathbf{0};\; \sum_{i=1}^n y_i = 1 \}$. Clearly, this optimization problem is concave in $\by$.

Finally, it is worth emphasizing that, the algorithms developed based on Danskin's theorem (or its variations) can only be applied to problems where step   (\ref{eq:ApproxgradientDescent}a) can be computed efficiently. While  non-convex non-concave min-max problems do not satisfy this assumption in general, one may be able to approximate the objective function with another objective function for which this assumption is satisfied. The following example illustrates this possibility in a particular application.

\medskip
\begin{example}\label{eq:exmaple:robust_training}
Consider the defense problem against adversarial attacks explained in the previous section where the training of a neural network requires solving the optimization problem (using similar notations as in \eqref{eq: poison_robust_attack}):
\begin{align}
\min_{\bw }\, \; \sum_{i=1}^{N}\;\max_{\|\boldsymbol{\delta}_i\|\leq \varepsilon}   {\ell}(p(\bu_i+\boldsymbol{\delta}_i;\bw ), t_i).\label{eq:robust-noncav}
\end{align}

 Clearly, the objective function is non-convex in $\mathbf{w}$ and non-concave in $\boldsymbol{\delta}$. Although finding the strongest adversarial attacker $\boldsymbol{\delta}_i$, that maximizes the inner problem in \eqref{eq:robust-noncav}, might be intractable,  it is usually possible to obtain a finite set of weak attackers. In practice, these attackers could be obtained using heuristics, e.g. projected gradient ascent or its variants~\cite{nouiehed2019solving}. Thus, \cite{nouiehed2019solving} proposes to approximate the above problem with the following more tractable version:
\begin{align} \label{eq:robust-cav}
\min_{\mathbf{w}}\; \frac{1}{N}\sum_{i=1}^{N}\;\max\Big\{\ell\left(p\left(\bu_{i}+\bdelta_k(\bu_{i},\bw);\mathbf{w}\right),t_i\right)\Big\}_{k=1}^{K},
\end{align}
where $\{\bdelta_k(\bu_i,\bw)\}_{k=1}^K$ is a set of $K$-weak attackers to data point $\bu_i$ using the neural network's weight $\bw$. Now the maximization is over a finite number of adversaries and hence can be transformed to a concave inner maximization problem using the transformation described in~\eqref{eq:finiteMax} and \eqref{eq:finiteMax2}. The performance  of this simple reformulation of the problem is depicted in  Table~\ref{tab:robust_nn_results}. As can be seen in this table, the proposed reformulation yields comparable results against state-of-the-art algorithms. In addition, unlike the other two algorithms, the proposed method enjoys theoretical convergence guarantees. \hfill $\blacksquare$

In the optimization problem~\eqref{eq:robust-noncav},  the inner optimization problem is non-concave in~$\boldsymbol{\delta}$. We approximate this non-concave function with a concave function by generating a finite set of adversarial instances in \eqref{eq:robust-cav} and used the transformation described in \eqref{eq:finiteMax},\eqref{eq:finiteMax2} to obtain a concave inner maximization probelm. This technique can be useful in solving other general  optimization problems of the form~\eqref{eq:MinMax} by approximating the inner problem with a finite set of points in the set $\mathcal{Y}$. Another commonly used technique to approximate the inner maximization in \eqref{eq:MinMax} with a concave problem  is to add a proper regularizer in $y$.  This technique has been used in \cite{sanjabi2018convergence} to obtain a stable training procedure for generative adversarial networks.



\end{example}
\begin{table}[]
\centering
\resizebox{0.5\textwidth}{!}{%
\begin{tabular}{@{}ccccc@{}}
\toprule
                                                       &                                        & A~\cite{madry2018}             & B~\cite{zhang_icml_2019}       & Proposed~\cite{nouiehed2019solving} \\ \midrule
\multicolumn{2}{c|}{Natural}                                                                    & \multicolumn{1}{c|}{98.58\%} & \multicolumn{1}{c|}{97.21\%} & 98.20\%  \\ \midrule
\multicolumn{1}{c|}{\multirow{3}{*}{FGSM~\cite{FGSM}}} & \multicolumn{1}{c|}{$\varepsilon=0.2$} & \multicolumn{1}{c|}{96.09\%} & \multicolumn{1}{c|}{96.19\%} & 97.04\%  \\
\multicolumn{1}{c|}{}                                  & \multicolumn{1}{c|}{$\varepsilon=0.3$} & \multicolumn{1}{c|}{94.82\%} & \multicolumn{1}{c|}{96.17\%} & 96.66\%  \\
\multicolumn{1}{c|}{}                                  & \multicolumn{1}{c|}{$\varepsilon=0.4$} & \multicolumn{1}{c|}{89.84\%} & \multicolumn{1}{c|}{96.14\%} & 96.23\%  \\ \midrule
\multicolumn{1}{c|}{\multirow{3}{*}{PGD~\cite{PGD}}}   & \multicolumn{1}{c|}{$\varepsilon=0.2$} & \multicolumn{1}{c|}{94.64\%} & \multicolumn{1}{c|}{95.01\%} & 96.00\%  \\
\multicolumn{1}{c|}{}                                  & \multicolumn{1}{c|}{$\varepsilon=0.3$} & \multicolumn{1}{c|}{91.41\%} & \multicolumn{1}{c|}{94.36\%} & 95.17\%  \\
\multicolumn{1}{c|}{}                                  & \multicolumn{1}{c|}{$\varepsilon=0.4$} & \multicolumn{1}{c|}{78.67\%} & \multicolumn{1}{c|}{94.11\%} & 94.22\%  \\ \bottomrule
\end{tabular}%
}
\caption{\footnotesize The performance of different defense algorithms for  training neural network on MNIST dataset. The first row is the accuracy when no attack is present. 
The second and the third row show the performance of different defense algorithms under ``FGSM" and ``PGD" attack, respectively.
Different $\varepsilon$ values show the magnitude of the attack. Different columns show different defense strategies. The defense method A (proposed in \cite{madry2018}) and the defense mechanism B (proposed in \cite{zhang_icml_2019}) are compared against the proposed method in~\cite{nouiehed2019solving}.
More details on the experiment can be found in~\cite{nouiehed2019solving}.}
\label{tab:robust_nn_results}
\end{table}
 
\subsection{Algorithms Based on Solving VI}
Another perspective that leads to the development of algorithms is the game theoretic perspective. To  present the ideas, first notice that the \eqref{eq:Game-Stat} notion defined in the previous section can be summarized as 
\begin{equation}\label{eq:VI}
\langle F(\bar{\bz}) , \bz - \bar{\bz}\rangle \geq 0, \;\forall \bz \in \mathcal{Z}
 \end{equation}
where $\mathcal{Z} \triangleq \mathcal{X} \times \mathcal{Y}$, 
$\bz = \left[
\begin{array}{c}
\bx\\
\by\\
\end{array}
\right]
$,
 and $F(\cdot)$ is a mapping induced by the objective function $f(\bx,\by)$ in \eqref{eq:MinMax}, defined by $$F(\bar{\bz}) \triangleq F\left(
  \left[
\begin{array}{c}
\bar{\bx}\\
\bar{\by}\\
\end{array}
\right]
 \right)=
  \left[
\begin{array}{r}
\nabla_{\bx}f(\bar{\bx},\bar{\by})\\
-\nabla_{\by}f(\bar{\bx},\bar{\by})\\
\end{array}
\right].$$ 
This way of looking at the min-max optimization problem naturally leads to the design of algorithms that compute the solution of the (Stampacchia) VI~\eqref{eq:VI}. 

When problem \eqref{eq:MinMax} is (strongly) convex in $\bx$ and (strongly) concave in $\by$, the mapping $F(\bz)$  is (strongly) monotone\footnote{A strongly monotone mapping $F(\cdot)$ satisfies the following $\langle F(\bz)-F(\bv),\bz-\bv\rangle\ge \sigma \|\bv-\bz\|^2,\; \forall~\bv,\bz\in \mathcal{Z}$, for some constant $\sigma>0$. If it satisfies this inequality for $\sigma = 0$, we say the VI is monotone.}, therefore classical methods for solving variational inequalities (VI)  such as extra-gradient can be applied \cite{nemirovski2004prox}. However, in the non-convex and/or non-concave setting of interest to this article, the strong monotonicity property no longer holds; and hence the classical algorithms cannot be used in this setting.  To overcome this barrier, a natural approach is to approximate the mapping $F(\cdot)$ with a series of strongly monotone mappings and solve a series of strongly monotone VIs. The work~\cite{lin2018solving} builds upon this idea and creates a series of strongly monotone VIs using proximal point algorithm, and  proposes an iterative procedure named inexact proximal point (IPP) method, as given below: 
\begin{subequations}\label{eq:ProximalPointVI}
\begin{align}
&\textrm{Let }     F^{\gamma}_{\bz^r}(\bz)=F(\bz)+\gamma^{-1}(\bz-\bz^r)\\
&\textrm{Let } \bz^{r+1} \textrm{  be the (approx) solution of VI } F^{\gamma}_{\bz^r}(\cdot).
\end{align}
\end{subequations}

In  \eqref{eq:ProximalPointVI}, $\gamma>0$ is chosen to be small enough so that the mapping $F^{\gamma}_{\bz^r}(\bz)$ becomes strongly monotone (in $\bz$). The strongly monotone mapping $F^{\gamma}_{\bz^r}(\bz)$ can be solved using another iterative procedure such as extra gradient method, or the iterative procedure 
\begin{equation} \label{eq:IterationVI}
    \bz^{t+1}=\mathcal{P}_{\mathcal{Z}}(\bz^t-\beta F(\bz^t))
\end{equation}
where $\beta$ denotes the stepsize and $\mathcal{P}_{\mathcal{Z}}$ is the projection to the set $\mathcal{Z}$. 

Combining the dynamics in \eqref{eq:ProximalPointVI} with the iterative procedure in~\eqref{eq:IterationVI} leads to a natural double-loop algorithm. This double-loop algorithm is not always guaranteed to solve the VI in \eqref{eq:VI}. Instead, it has been shown that this double-loop procedure computes a solution $\bz^*$ to the following {\it Minty VI}: 
\begin{equation}\label{eq:minty}
    \langle F(\bz),\bz-\bz^*\rangle\ge 0, \; \forall\bz\in\mathcal{Z}.
\end{equation}
Notice that this solution concept is different than the solution $\bar\bz$ in \eqref{eq:VI} as it has $F(\bz)$ instead of $F(\bar\bz)$ in the left hand side.
While these two solution concepts are different in general, it is known that if $\mathcal{Z}$ is a convex set, then any $\bz^*$ satisfying \eqref{eq:minty} also satisfies \eqref{eq:VI}. Furthermore, if 
$F(\cdot)$ is monotone (or  when $f(\cdot,\cdot)$ is convex in $\bx$ and concave in $\by$), then any solution to  \eqref{eq:VI} is also a solution to \eqref{eq:minty}. 
While such monotonicity requirement can be slightly relaxed to cover a wider range of non-convex problems (see e.g. \cite{qihang18}), it is important to note that for generic non-convex, and/or non-convave  function $f(\cdot, \cdot)$,  there may not exist $\bz^*$ that satisfies \eqref{eq:minty}; see Example~\ref{eq:exmaple:1}. 

\medskip
\begin{example}\label{eq:exmaple:1}
Consider the following function which is non-convex in $x$, but concave in $y$:
$$f(x,y)=x^3+2xy-y^2, \; \mathcal{X}\times \mathcal{Y} = [-1,1]\times[-1,1].$$
One can verify that there are only two points, $(0,0)$ and $(-1,-1)$ that satisfy \eqref{eq:VI}. However, none of the above solutions satisfies \eqref{eq:minty}. To see this, one can verify that  $\langle F(z), z - z^* \rangle= -4 <0$ for $z=(0,-1)$ and $z^* = (-1,-1)$ and that $\langle F(z), z - z^* \rangle= -3 <0$ for $z=(-1,0)$ and $z^* = (0,0)$. Since any $\bz^*$ satisfying \eqref{eq:minty} will satisfy \eqref{eq:VI}, we conclude that there is no point satisfying \eqref{eq:minty} for this min-max problem. \hfill $\blacksquare$
\end{example}

In conclusion, the VIs \eqref{eq:VI} and \eqref{eq:minty} offer new perspectives to analyze problem~\eqref{eq:MinMax}; but the existing algorithms such as  \eqref{eq:ProximalPointVI} cannot deal with many problems covered by the potential based methods discussed in Sec.~\ref{sub:potential} (for example when $f(\bx,\by)$ is non-convex in $x$ and concave in $\by$ or when a point satisfying \eqref{eq:minty} does not exists as we explained in Example~\ref{eq:exmaple:1}). Moreover, the VI-based algorithms completely ignore the order of maximization and minimization in \eqref{eq:MinMax} and, hence, cannot be applied to problems in which the order of min and max is crucial.

\subsection{Algorithms Using Single-Loop Update}\label{eq:alternating}
The algorithms discussed in the previous two subsections are all {\it double loop} algorithms, in which one variable (e.g. $\by$) is updated in a few consecutive iterations before another variable gets updated. 
In many practical applications, however, {\it single loop} algorithms which update $\bx$ and $\by$ either alternatingly  or simultaneously are  preferred. For example, in problem \eqref{eq:jammer}, the jammer often pretends to be the regular user, so it updates simultaneously with the regular users \cite{gohary09}.
However, it is challenging to design and analyze single loop algorithms for problem \eqref{eq:MinMax} --- even for the simplest linear problem where $f(\bx,\by) = \langle \bx,\by\rangle$, the single-loop algorithm GDA diverges; see the discussion in Sec. \ref{sec:challenges}.

To overcome the above challenges,  \cite{minmaxtsp19}  proposes a single loop algorithm called Hybrid Block Successive Approximation (HiBSA), whose iterations are given by 
\begin{subequations}
\begin{align}
\hspace{-0.2cm}&\mathbf {x}^{r+1} \hspace{-0.1cm}  = \mathcal{P}_{\mathcal X} \hspace{-0.1cm} \left  ( \mathbf {x}^{r}- \beta^r {\nabla}_{\mathbf x} f \left (\mathbf {x}^{r},\mathbf y^{r} \right ) \right ),  
\\
\hspace{-0.2cm}&\mathbf {y}^{r+1} \hspace{-0.1cm}  = \mathcal{P}_{\mathcal Y} \hspace{-0.1cm} \left  ( \mathbf (1+\gamma^r\rho){\by}^{r} \hspace{-0.1cm} + \rho {\nabla}_{\mathbf y} f \left (\mathbf {x}^{r+1},\mathbf y^{r} \right ) \right ), \label{eq.updatofy}
\end{align}
\end{subequations}
where $\beta^r, \rho > 0$ are the step sizes; $\gamma^r>0$ is some perturbation parameter. This algorithm can be further generalized to optimize certain approximation functions of $x$ and $y$, similarly to the successive convex approximation strategies used in min-only problems \cite{meisam15BSUMsurvey,Razaviyayn12SUM}. The ``hybrid" in the name refers to the fact that this algorithm contains both the descent and ascent steps.

The HiBSA iteration is very similar to the GDA algorithm mentioned previously, in which gradient descent and ascent steps are performed alternatingly. The key difference here is that the $\by$ update includes an additional term $\gamma^r\rho\by^{r}$, so that at each iteration the $\by$ update represents a ``perturbed" version of the original gradient ascent step. The idea is that after the perturbation, the new iteration $\by^{r+1}$ is ``closer" to the old iteration $\by^{r}$, so it can avoid the divergent patterns depicted in Fig.~\ref{Nowazin17}. 
Intuitively, as long as the perturbation eventually goes to zero, the algorithm will still converge to the desired solutions. Specifically,
it is shown in \cite{minmaxtsp19} that, if $f(\bx,\by)$ is strongly concave in $\by$, then one can simply remove the perturbation term (by setting $\gamma^r=0$ for all $r$), and the HiBSA will converge to a point satisfying  condition~\eqref{eq:Game-Stat}. Further, if $f(\bx,\by)$ is only concave in $\by$, then one needs to choose  $\beta^r = \mathcal{O}(1/{r}^{1/2})$, and $\gamma^r=\mathcal{O}(1/r^{1/4})$ to converge to a point satisfying   condition~\eqref{eq:Game-Stat}. 

\begin{figure}[htb]
	\begin{center}
	\begin{minipage}[t]{0.30\textwidth}
	\includegraphics[width=\textwidth]{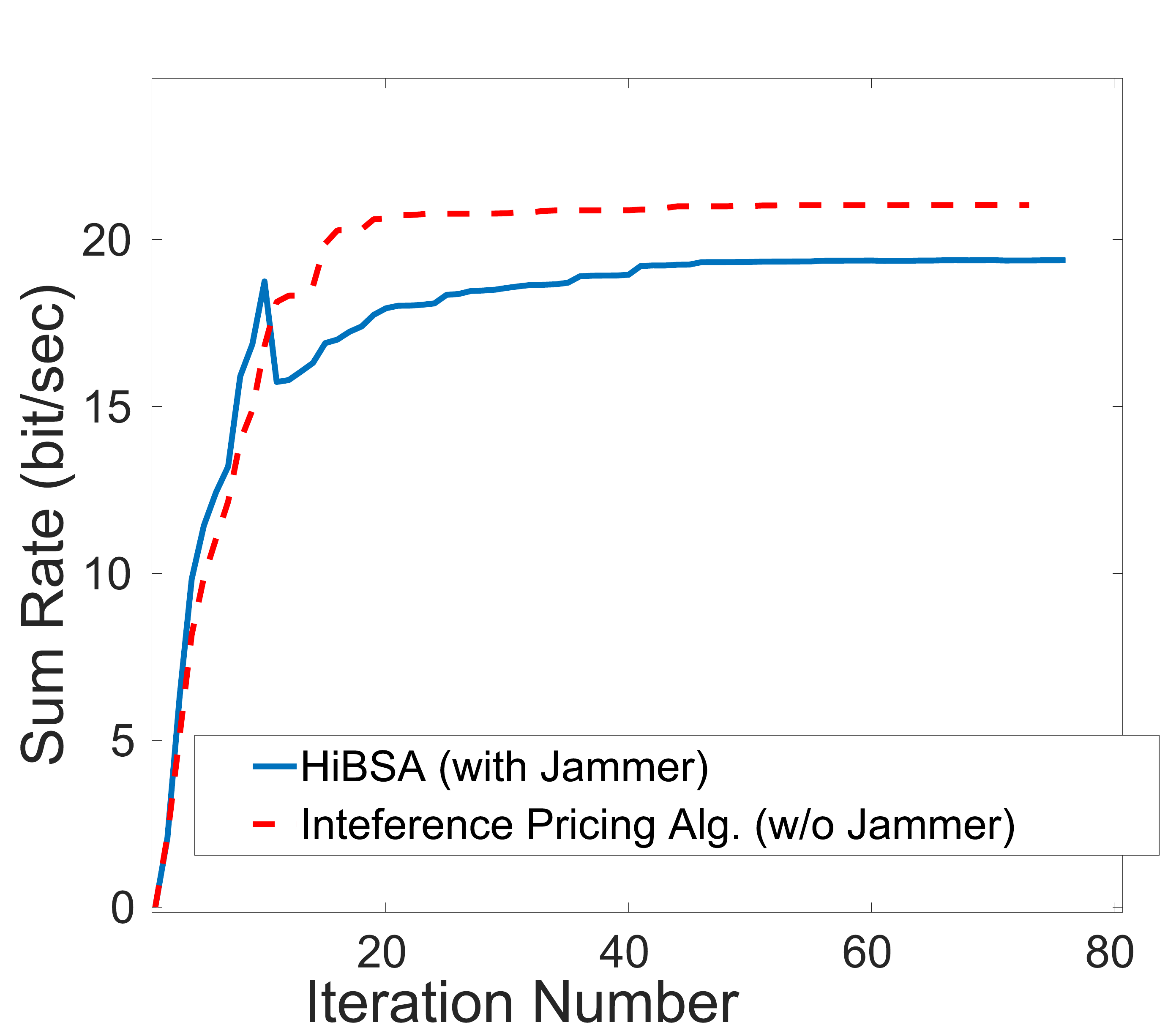}
		\end{minipage}
	\begin{minipage}[t]{0.345\textwidth}
	\includegraphics[width=\textwidth]{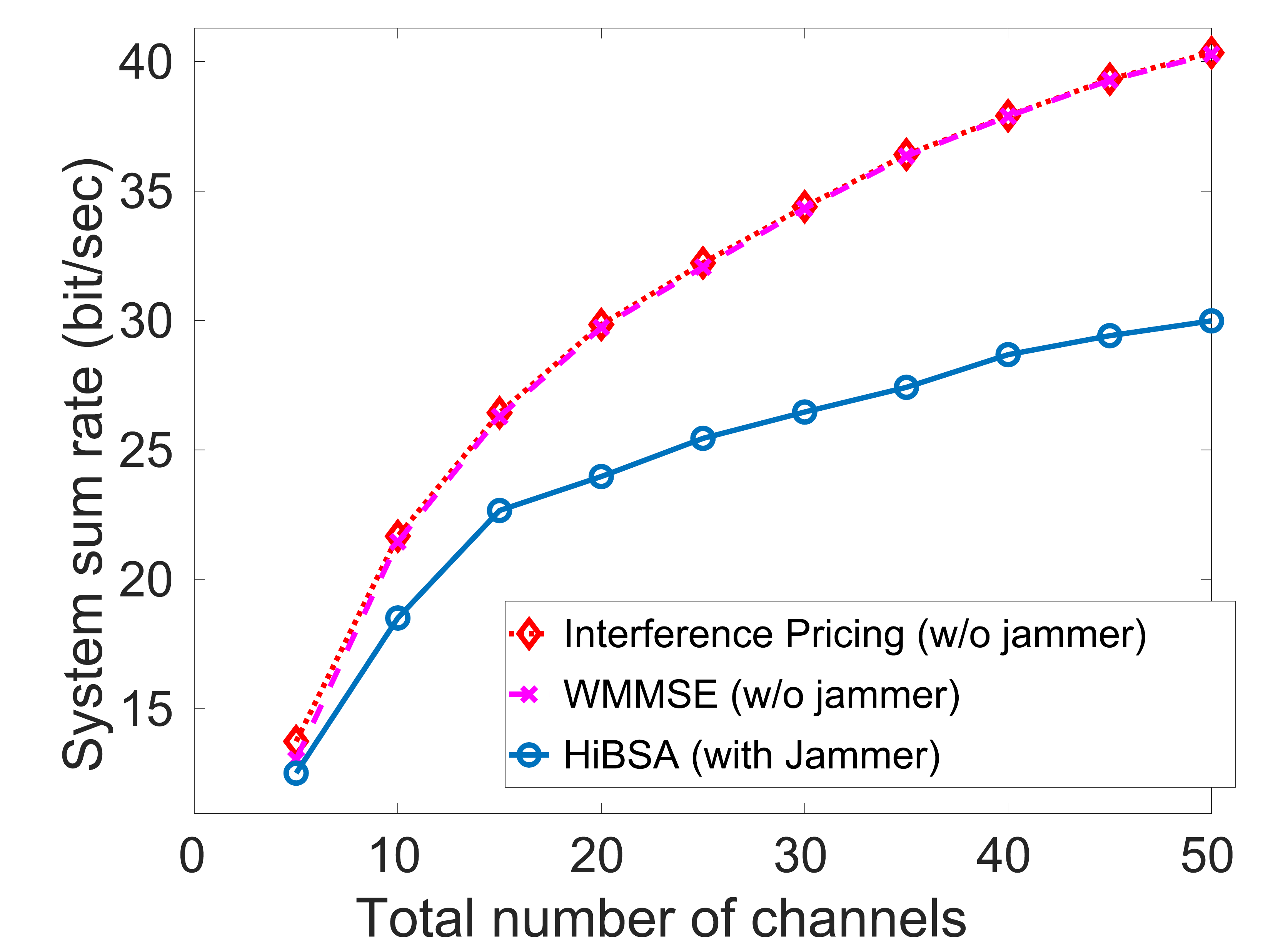}
		\end{minipage}
	\end{center}
	\caption{\footnotesize The convergence curves and total averaged system performance comparing three algorithms: WMMSE, Interference Pricing and HiBSA. All  users' power budget is fixed at $P=10^{\rm SNR/10}$. For test cases without a jammer, we set $\sigma^2_k=1$ for all $k$. For test cases with a jammer, we set $\sigma^2_k=1/2$ for all $k$ and let the jammer have the rest of the noise power, i.e., $p_{0,\max}= N/2$.  
	Figure taken from \cite{minmaxtsp19}.} 
	\label{fig:jamming}
\end{figure}

\medskip
\begin{example}
We apply the HiBSA to the power control problem \eqref{eq:jammer}. 
It is easy to verify that the jammer's objective is strongly concave over the feasible set.

We compare HiBSA with two  classic algorithms: interference pricing \cite{Schmidt13compare}, and the WMMSE \cite{shi11WMMSE_TSP}, both of which are designed for power control problem {\it without} assuming the presence of the jammer. Our problem is tested using the following setting. We construct a network with $10$ regular user and a single jammer. The interference channel among the users and the jammer is generated using  uncorrelated fading channel model with channel coefficients generated from the complex zero-mean Gaussian distribution with unit covariance.  

From Fig.~\ref{fig:jamming} (top), we see that the pricing algorithm monotonically increases the sum rate (as is predicted by theory), while HiBSA behaves differently: after some initial oscillation, the algorithm converges to a value that has lower sum-rate. Further in Fig.~\ref{fig:jamming} (bottom), we do see that by using the proposed algorithm, the jammer is able to effectively reduce the total sum rate of the system. \hfill $\blacksquare$
\end{example}



\subsection{Extension to Zeroth-order Based Algorithms}
Up to now, all the algorithms reviewed require first-order  (gradient) information. In this subsection, we disucss a useful extension when only {\it zeroth-order} (ZO) information is available. That is, we only have access to the objective values $f(\bx,\by)$ at a given point $(\bx,\by)$ at every iteration. This type of algorithm is useful, for example, in practical adversarial attack scenario where the attacker only has access to the output of the ML model~\cite{zominmax19}. 


To design algorithms in the ZO setting, one typically replaces the gradient  $\nabla h(\bx)$ with some kind of {\it gradient estimate}. One popular estimate is given by 
$$\widehat{\nabla}_{\mathbf x}h  (\mathbf x ) =
    {   \frac{1}{q  } \sum_{i=1}^q  \frac{d [h ( \mathbf x + \mu \mathbf u_i) - h ( \mathbf x  ) ] }{\mu} }  \mathbf u_i,$$
where $\{ \mathbf u_i \}_{i=1}^q$ are $q$ \emph{i.i.d.} random direction vectors drawn uniformly from the unit sphere, and $\mu > 0$ is a  smoothing parameter.  We note that the ZO gradient estimator involves the random direction sampling w.r.t.  $\mathbf u_i$.
It is known that    $\widehat{\nabla}_{\mathbf x}h (\mathbf x ) $ provides an unbiased estimate of the gradient of the smoothing function of $f$ rather than 
the true gradient of $f$. Here   the smoothing function  of $f$ is defined by
$h_{\mu}(\mathbf x) = \mathbb E_{\mathbf v} [ h(\mathbf x + \mu \mathbf v) ]$,  where
 $\mathbf v$ follows the uniform distribution over the unit Euclidean ball. Such a gradient estimate is used in~\cite{zominmax19} to develop a ZO algorithm for solving min-max problems.

\begin{figure}[htb]
 \begin{center}
 \begin{minipage}[t]{0.32\textwidth}
 \includegraphics[width=\textwidth]{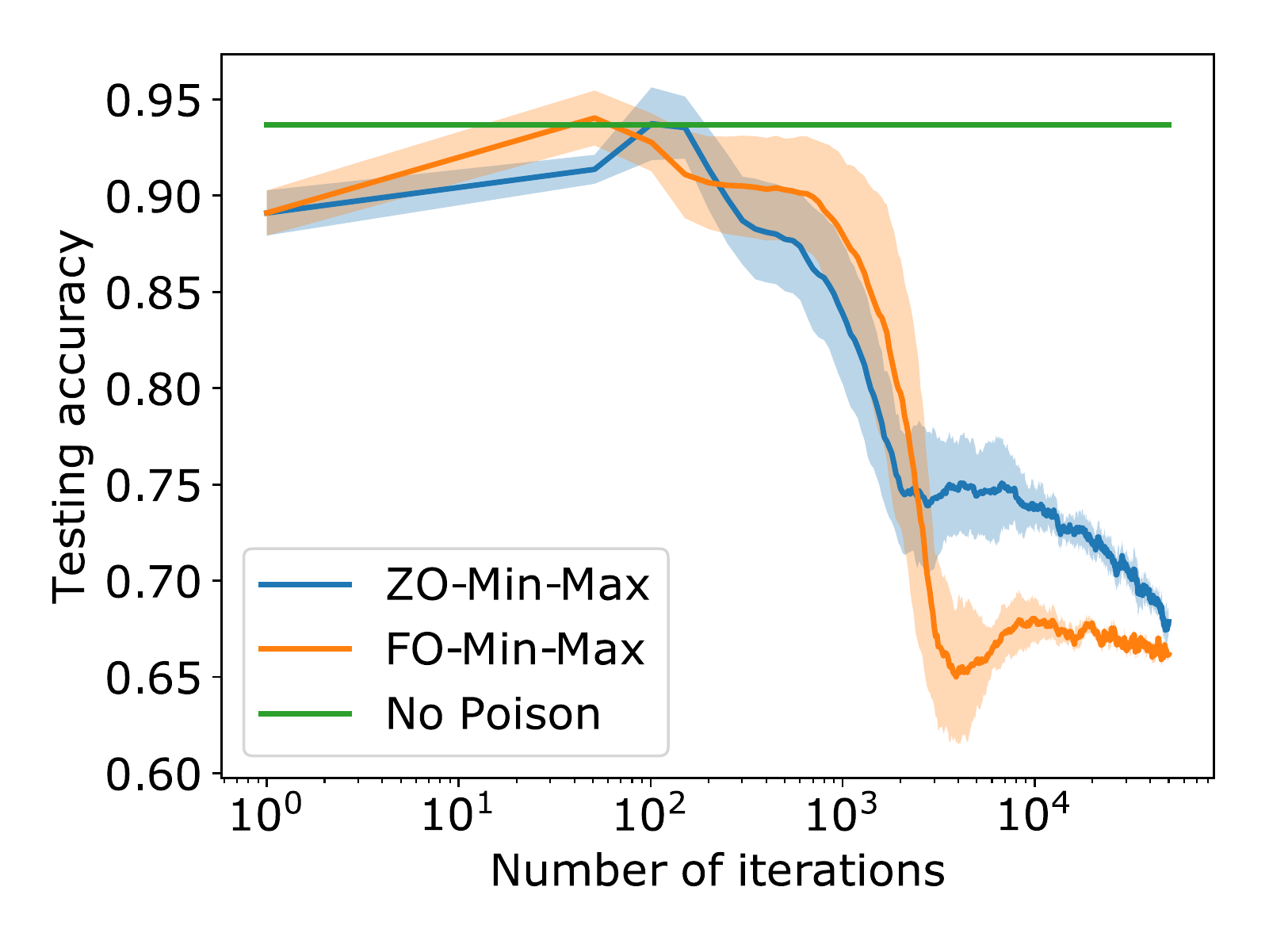}
 \end{minipage}
 \begin{minipage}[t]{0.32\textwidth}
 \includegraphics[width=\textwidth]{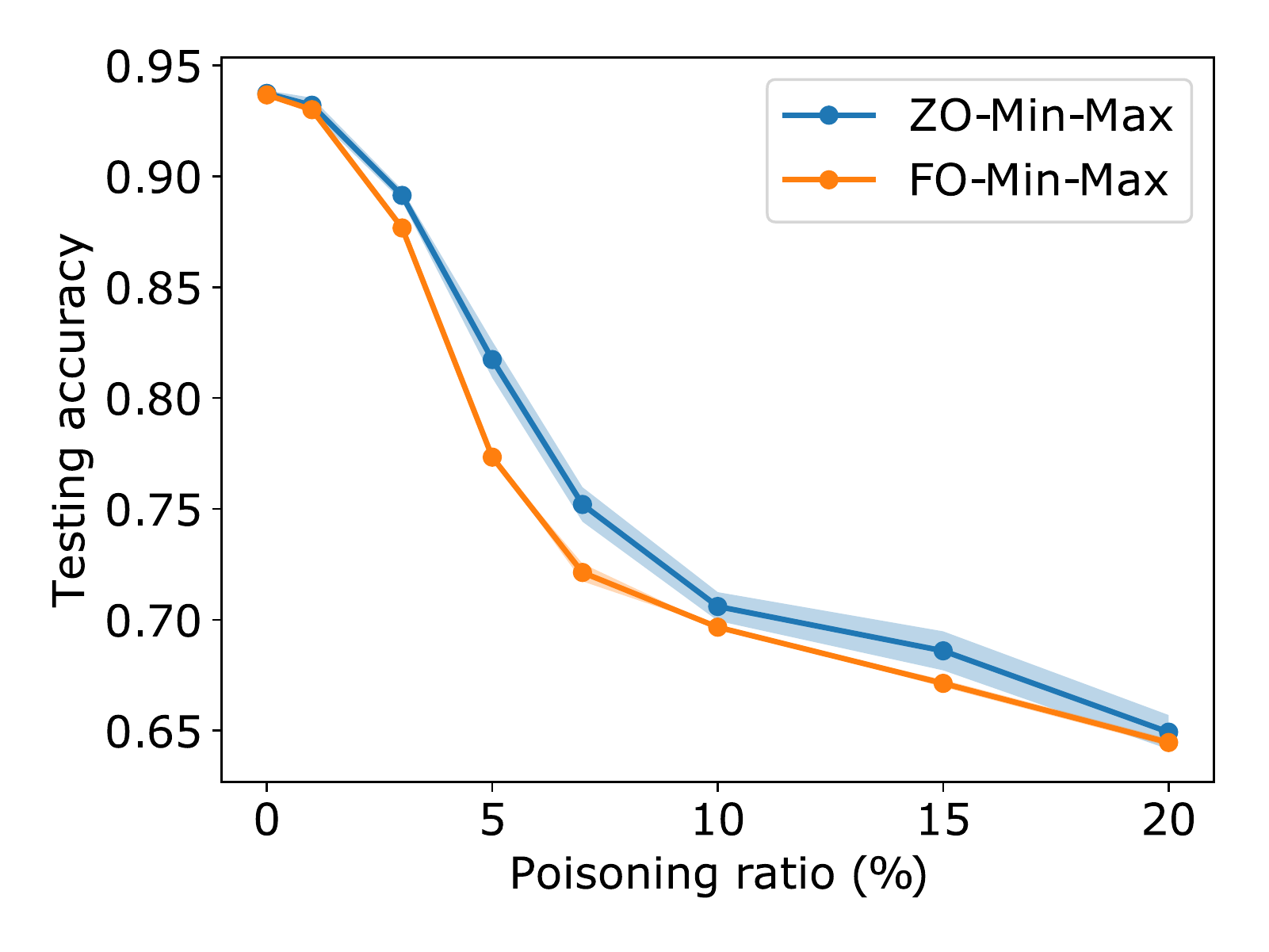}
 \end{minipage}
 \end{center}
 \caption{\footnotesize Empirical performance of ZO/FO-Min-Max in poisoning attacks, where the ZO/FO-Min-Max algorithm refers to HiBSA with either ZO or FO oracle: (top) testing accuracy  versus iterations (the shaded region represents variance of $10$ random trials), and
 (bottom) testing accuracy versus  data poisoning ratio. 	Figure taken from \cite{zominmax19}.} 
 \label{fig: poison_attack}
 \end{figure}

\medskip
\begin{example}
To showcase the performance comparison between ZO based and first-order (FO) based algorithm, we consider applying HiBSA and its ZO version to the data poisoning problem \eqref{eq: poison_robust_attack}.  
In particular, let us consider attacking the  data set used to train a logistic regression model.  
We first set the poisoning ratio, 
i.e.,  the percentage of the training samples attacked, to $15\%$.
Fig.\,\ref{fig: poison_attack} (top) demonstrates the testing accuracy (against iterations) of the  model learnt from   poisoned training  data, where the poisoning attack is generated by ZO min-max (where the adversarial only has access to victim model outputs) and FO min-max (where the adversarial has access to details of the victim model). 
As we can see from Fig. \ref{fig: poison_attack}, the poisoning attack can significantly reduce the testing accuracy compared to the clean model. Further, the ZO min-max yields promising attacking performance comparable to the FO min-max. 
Additionally, in Fig.\,\ref{fig: poison_attack} (bottom), we present the testing accuracy of the learned model under different data poisoning ratios. As we can see,  only $5\%$ poisoned training data can significantly break the testing accuracy of a well-trained model. The details of this experiment can be found in \cite{zominmax19}. \hfill $\blacksquare$
\end{example}

\section{Connections among algorithms and optimality conditions} 
\begin{figure}[htb]
	\begin{center}
	\includegraphics[width=0.5\textwidth]{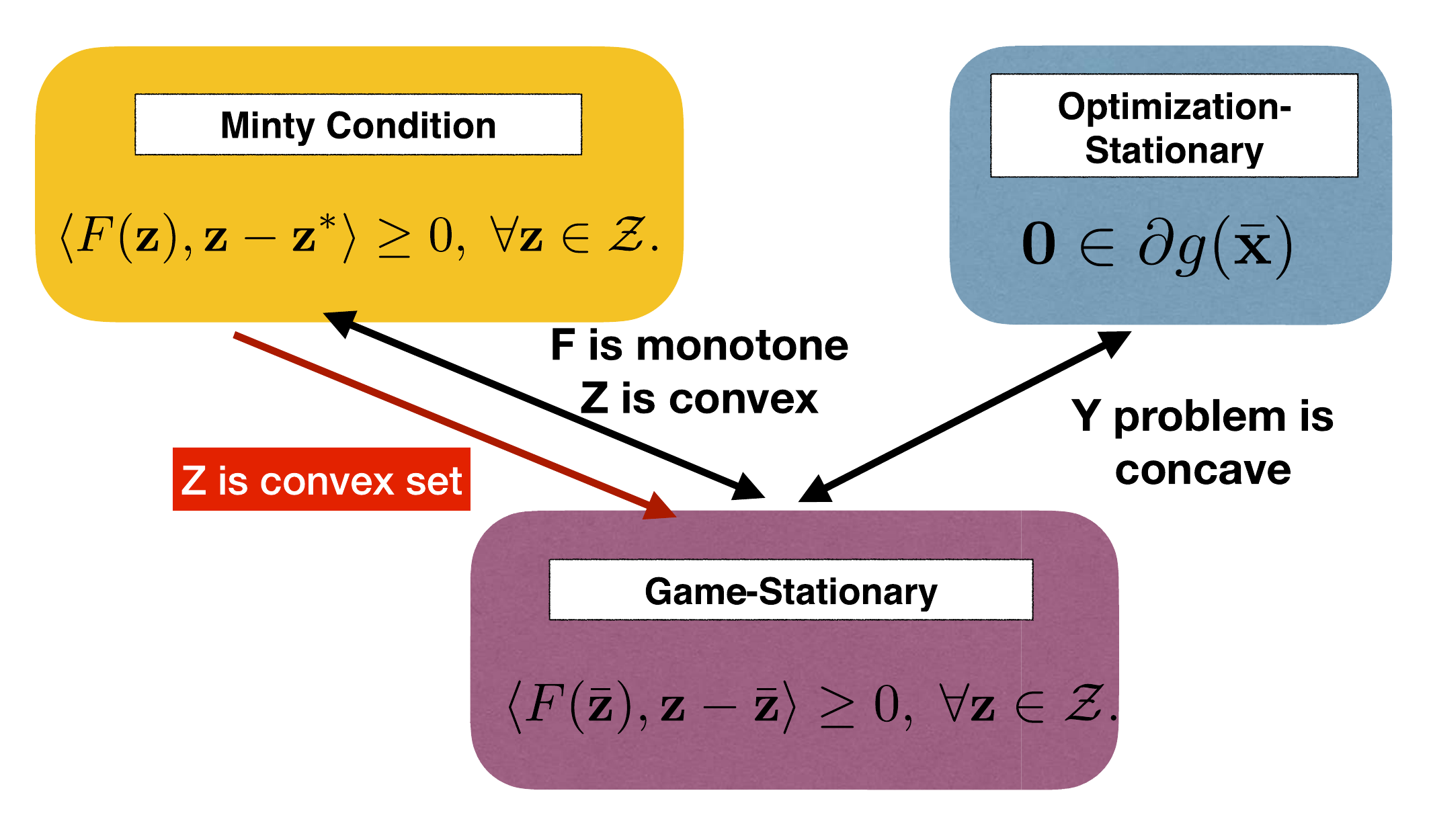}
	\end{center}
	\caption{\footnotesize Relations of different optimality conditions}. 
	\label{fig:optimality}
\end{figure}

{In this section, we summarize our discussion on various optimality conditions as well as algorithm performance. First, in Fig. \ref{fig:optimality}, we describe the relationship between the Minty condition \eqref{eq:minty}, the \eqref{eq:Opt-Stationary} and \eqref{eq:Game-Stat}. 
 Second, we compare the properties of different algorithms, such as their convergence conditions and optimality criteria in Table~\ref{tb:rate}. Despite the possible equivalence between the optimality conditions, we still keep the column ``optimality criteria" because these are the criteria based on which the algorithms are originally designed.} 

\begin{table*}
{\footnotesize
\begin{center}
\begin{threeparttable}
\begin{tabular}{|c|c|c|c|c|c|c|}
\hline
\textbf{Algorithm} & \textbf{ Optimality Criterion} & \textbf{Oracles} & \multicolumn{2}{c|}{\textbf{Assumptions}} & \textbf{Other Comments}   \\
\hline
\multirow{2}{*}{ \multirow{2}{*}\makecell{ \textit{Multi-Step GDA}  \cite{nouiehed2019solving}} }  &  \multirow{2}{*}{ \makecell { Game-Stationary}} &  \multirow{2}{*}{ \makecell { FO }}  & \multicolumn{2}{c|}{\multirow{2}{*}{\makecell{ $f$ NC in $x$, concave in $y$} }} & {\multirow{2}{*}{\makecell{  Det. \& DL, \& Acc.}}}  \\
&  & & \multicolumn{2}{c|}{}  & 
\\ \hline
\multirow{2}{*}{ \multirow{2}{*}\makecell{ \textit{CDIAG}  \cite{NIPS2019_9430}} }  &  \multirow{2}{*}{ \makecell { Optimization-Stationary}} &  \multirow{2}{*}{ \makecell { FO }}  & \multicolumn{2}{c|}{\multirow{2}{*}{\makecell{ $f$ NC in $x$, concave in $y$} }} & {\multirow{2}{*}{\makecell{  Det. \& TL, \& Acc.}}}  \\
&  & & \multicolumn{2}{c|}{}  & 
\\ \hline

\multirow{4}{*}{\makecell{ \\ \\ \textit{PG-SMD/ PGSVRG}\cite{qihang18} \\}}    & \multirow{4}{*}{ \makecell { \\ \\ Optimization-Stationary\\}} & \multirow{4}{*}{ \makecell { \\ \\ FO \\  }} &  \multirow{2}{*}{\makecell } & \makecell{\\  $f$ concave in $y$} & \multirow{4}{*}{ \makecell { \\ \\ St. \& DL\\}}    \\ 
& & & &  & 
\\ \cline{4-5}
& & & \multirow{2}{*}{\makecell{$f = \frac{1}{n} \sum_{i=1}^{n} f_{i}$ \\[1ex]  $f_i(x,y)$  NC in $x$}} & \makecell{\\  $f_{i}$ concave in $y$}  & \\ 
& & & &  & 
\\ \hline

\multirow{2}{*}{ \makecell{ \textit{IPP} for VI \cite{lin2018solving}} }  &  \multirow{2}{*}{\makecell {Minty-Condition}} &  {\multirow{2}{*}{\makecell{ FO}}}     & \multicolumn{2}{c|}{\multirow{2}{*}{\makecell{  NC in x, y}}}  & {\multirow{2}{*}{\makecell{  Det. \& DL}}}  
\\ 
&  & & \multicolumn{2}{c|}{}  & \\ \hline

\multirow{2}{*}{ \makecell{ \textit{GDA }\cite{lin2019gradient}} }  &  \multirow{2}{*}{\makecell { Optimization-Stationary}} &  {\multirow{2}{*}{\makecell{ FO}}}     & \multicolumn{2}{c|}{\multirow{2}{*}{\makecell{ $f$ NC in $x$, concave in $y$}}}  & {{\makecell{  Det. \& St. \& DL}}}  
\\ 
&  & & \multicolumn{2}{c|}{}  & {$y$-step has small stepsize}\\ \hline

\multirow{2}{*}{ \makecell{ \textit{HiBSA }\cite{minmaxtsp19}\cite{zominmax19}} }  &  \multirow{2}{*}{\makecell { Game-Stationary}} &  {\multirow{2}{*}{\makecell{ FO \& ZO}}}     & \multicolumn{2}{c|}{\multirow{2}{*}{\makecell{ $f$ NC in $x$, concave in $y$}}}  & {\multirow{2}{*}{\makecell{  Det. \& SL}}}  
\\ 
&  & & \multicolumn{2}{c|}{}  & \\ \hline

\end{tabular}
\end{threeparttable}
\end{center}
\caption{{\footnotesize {Summary of algorithms for the min-max optimization problem \eqref{eq:MinMax} along with their convergence guarantees. Note that in the third column, we characterize the type of the oracle used, i.e., FO or ZO. In the last column are other comments about the algorithms}, i.e deterministic (Det.) or stochastic (St.), single loop (SL) or double loop (DL) or triple loop (TL), Acceleration (Acc.) or not. Moreover, we use the abbreviations NC for non-convex. }}
\label{tb:rate}
}
\end{table*}


\section{Conclusion and Future research directions}
Non-convex min-max optimization problems appear in a wide range of applications. 
Despite the recent developments in solving these problems, the available tool sets and theories are still very limited. In particular, as discussed in section~\ref{sec:Algorithms}, these algorithms require at least one of the two following assumptions:
\begin{itemize}
    \item[i)] The objective of one of the two player is easy to optimize. For example, the object function in \eqref{eq:MinMax} is concave in~$\by$.
    \item[ii)] The Minty solutions satisfying~\eqref{eq:minty} for the min-max game are the solutions to the Stampchia VI~\eqref{eq:VI}.
    \end{itemize}
While the first assumption is mostly easy to check, the second assumption might not be easy to verify.  Nevertheless, these two conditions do not imply each other. Moreover, there is a wide range of non-convex min-max problem instances that does not satisfy either of these assumptions. For solving those problems, it might be helpful to approximate them with a min-max problem satisfying one of these two assumptions. 

As future research directions, a natural first open direction is toward the development of algorithms that can work under a more relaxed set of assumptions. We emphasize that the class of problems that are provably solvable (to satisfy either \eqref{eq:Opt-Stationary} or \eqref{eq:Game-Stat}) is still very limited, so it is important to extend solvable problems to a more general set of non-convex non-concave functions. One possible first step to address this is to start from algorithms that converge to desired solutions when initialized close enough to them, i.e. {\it local} convergence; for recent developments on this topic see \cite{adolphs2019local}. Another natural research direction is about the development of the algorithms in the absence of smoothness assumptions. When the objective function of the players are non-smooth but ``proximal-gradient friendly'', many of the results presented in this review can still be used by simply using proximal gradient instead of gradient~\cite{barazandeh2020solving}. These scenarios happen, for example, when the objective function of the players is a summation of a smooth function and a convex non-smooth function.  Additionally, it is of interest to customize the existing generic non-convex min-max algorithms to practical applications, for example to reinforcement learning  \cite{NIPS2019_8814}.

One of the main applications of non-convex min-max problems, as mentioned in section~\ref{sec:applications},  is to design  systems that are robust against uncertain parameters or the existence of adversaries. A major question in these applications is whether one can provide ``robustness certificate" after solving the optimization problem. In particular, can we guarantee or measure the robustness level in these applications? The answer to this question is closely tied to the development of algorithms for solving non-convex min-max problems. 

Another natural research direction is about the rate of convergence of the developed algorithms. For example, while we know solving min-max problems to \eqref{eq:Opt-Stationary} is easy when \eqref{eq:MinMax} is concave in $\by$ (and possibly non-convex in $\bx$), the optimal rate of convergence (using gradient information) is still not known. Moreover, it is natural to ask whether knowing the Hessian or higher order derivatives of the objective function could improve the performance of these algorithms. So far, most of the algorithms developed for non-convex min-max problems rely only on gradient information.


\bibliographystyle{IEEEtran}
\bibliography{ref}

 \clearpage 

\end{document}